\documentclass[11pt,a4paper]{article}
\usepackage{fullpage}
\usepackage{authblk}

\usepackage[figuresright]{rotating}
\usepackage{amssymb}
\usepackage{amsthm}
\usepackage{multicol}
\usepackage{subfigure}
\usepackage{graphicx}
\usepackage{epstopdf}
\usepackage{fullpage}
\usepackage{latexsym,amsmath}
\usepackage{stmaryrd,appendix}
\usepackage{algorithm,algorithmic}
\usepackage{subfigure}
\usepackage{caption}
\usepackage{amsfonts}
\usepackage{xcolor,multirow}
\usepackage{cite}
\usepackage{bm}
\usepackage{url}
\usepackage{diagbox}
\usepackage{array}
\usepackage{booktabs}

\newtheorem{lemma}{Lemma}
\newtheorem{theorem}{Theorem}

\title{Color image recovery using low-rank quaternion matrix completion algorithm}
\author{Jifei Miao\thanks{jifmiao@163.com} }

\author{Kit Ian Kou\thanks{kikou@umac.mo}}

\affil{\normalsize{Department of Mathematics, Faculty of Science
		and Technology, University of Macau, Macau 999078, China}}

\date{}

\begin{document}
	\maketitle
	\begin{abstract}
		\normalsize
 As a new color image representation tool, quaternion has achieved excellent results in color image processing problems. In this paper, we propose a novel low-rank quaternion matrix completion algorithm to recover missing data of color image. Motivated by two kinds of low-rank approximation approaches (low-rank decomposition and nuclear norm minimization) in traditional matrix-based methods, we combine the two approaches in our quaternion matrix-based model. Furthermore, the nuclear norm of the quaternion matrix is replaced by the sum of Frobenius norm of its two low-rank factor quaternion matrices. Based on the relationship between quaternion matrix and its  equivalent complex matrix, the problem eventually is converted from quaternion number field to complex number field. An alternating minimization method is applied to solve the model. Simulation results on real world color image recovery show the superior performance and efficiency of the proposed algorithm over some state-of-the-art tensor-based ones.
\end{abstract}


\section{Introduction}

Color images are applied in numerous fields, from the casual documentation of events to medical applications. A color image contains red, blue, and green channels. In most cases, some data of the acquired images are missed during acquisition and transmission.
Hence, a well-performed recovery technology should be proposed.

In the past few decades, low-rank matrix completion problem has been widely studied and proven very useful in the application of image recovery \cite{DBLP:journals/tit/KeshavanMO10, DBLP:journals/pami/HuZYLH13, DBLP:journals/ijon/WangZC14, DBLP:journals/sigpro/LinW15, DBLP:journals/mssp/YuPY19}.
Commonly, the method is to stack all the image pixels as column vectors of a matrix, and recovery theory and algorithm are adopted to the resulting matrix which is low-rank or  approximately low-rank.
However, these image recovery models are usually developed for gray-level images. For color image processing,  traditional matrix-based methods usually ignore the mutual connection amoung channels, because these recovery methods are applied to red, green, and blue channels separately, which is likely to result in color distortion during the recovery process.

More recently, an increasing number of low-rank tensor completion methods have been proposed to recover color images \cite{DBLP:journals/tsp/YokotaZC16, DBLP:journals/pami/ZhaoZC15, DBLP:journals/pami/LiuMWY13,  DBLP:journals/tip/BenguaPTD17, DBLP:journals/tip/ZhouLLZ18,DBLP:journals/sigpro/LongLCZ19}. Actually, a tensor is a natural form of high-dimensional and multi-way real world data. For example, a color image can be regarded as a 3-way tensor due to its three channels, each frontal slice of this 3-way tensor corresponds to a channel of the color image. The state-of-the-art methods for tensor completion mainly consist of approaches
of two types. The first type is nuclear norm minimization \cite{DBLP:journals/pami/LiuMWY13, DBLP:journals/tsp/ZhangA17, DBLP:conf/cvpr/ZhangEAHK14}, which is generally computationally expensive and time-consuming, for example some of this type of algorithms require computing tensor singular value decomposition (t-SVD) which costs much computation especially for natural large scale data. The second type involves the use of low-rank tensor decomposition techniques, for example some  Tucker decomposition based techniques have been proposed in \cite{DBLP:journals/pami/ChenHL14, DBLP:journals/pami/LiuMWY13}, and some CP (CANDECOMP/PARAFAC) decomposition based techniques have been proposed in \cite{DBLP:journals/pami/ZhaoZC15, DBLP:journals/tsp/YokotaZC16}. Nevertheless, for this type of algorithms, the rank of a tensor is generally pretty hard to determine \cite{DBLP:journals/tip/ZhouLLZ18}, so they usually cannot offer the best low-rank approximation to a tensor.
In brief, the recovery theory for low-rank tensor completion problem is not well established compared with that of matrix-based completion problem.

Different from conventional matrix-based and tensor-based models, in this paper, we propose a novel low-rank quaternion matrix completion algorithm to recover missing data of color image. Actually, the Red, Green and Blue values of each pixel of a color image can be naturally represented as a single pure quaternion valued pixel \cite{541817}.

\begin{equation}
\label{equ1}
q(x,y)=r(x,y)i+g(x,y)j+b(x,y)k,
\end{equation}
where $r(x,y)$, $g(x,y)$ and $b(x,y)$ are, respectively, the red,
green and blue components corresponding to the pixel at
position $(x,y)$ in the color image, and $i$, $j$ and $k$ are the three
imaginary units. By using (\ref{equ1}), an $M\times N$ color
image is described by a matrix with size $M\times N$ whose elements are pure quaternions.
The main advantage of this representation is that it processes a colour image holistically as a vector field and handles the coupling between the color channels naturally \cite{DBLP:journals/iet-ipr/ChenLSLS14, DBLP:journals/ijcv/SubakanV11, DBLP:journals/sigpro/ChenSZCTDL12}, and color information of source image is fully used. Moreover, comparing to the tensor-based model, the quaternion-based model not only preserves the correlation among channels but also the orthogonal property for the coefficients of different channels, which achieves a structured representation \cite{DBLP:journals/tip/XuYXZN15}. Hence, as a new color image representation tool, quation has achieved excellent results in the color image processing including color image filtering \cite{DBLP:journals/iet-ipr/ChenLSLS14}, color image edge detection \cite{DBLP:journals/imst/XuYL10}, color image denoising \cite{DBLP:journals/mssp/GaiYW015},  color image watermarking \cite{1199526_P2003}, color face recognition \cite{DBLP:journals/tip/ZouKW16}, color image recovery \cite{6744018} and so on.

However, the quaternion matrix completion based color image recovery problem  has  been  less investigated. In \cite{6744018}, the authors proposed a quaternion matrix completion algorithm by solving a semi-definite programming optimization (SDP) problem which can be solved by the interior-point method. In many practical applications, nevertheless, the matrices are very large, which makes the SDP problem intractable \cite{DBLP:journals/pami/HuZYLH13}. When the size of the matrices exceeds $100\times 100$, the existing state-of-the-art SDP solvers such as SDPT3 \cite{Toharticle2012} and SeDuMi \cite{Sturmarticle1998} are generally no longer applicable.
 In this paper, based on low-rank decomposition of quaternion matrix, we propose a novel low-rank matrix completion algorithm in quaternion number field. To the best of our knowledge, the quaternion matrix completion problem based on low-rank decomposition has not yet been investigated.  In summary, our main contributions include:
\begin{itemize}
	\item We propose an efficient low-rank quaternion matrix completion algorithm to recover missing data of color image. Both low-rank decomposition and nuclear norm minimization techniques are combined in our quaternion matrix-based model. The nuclear norm of the quaternion matrix is replaced by the sum of Frobenius norm of its two low-rank factor quaternion matrices. Furthermore, based on the relationship betwween the quaternion matrix and its  equivalent complex matrix, the problem eventually be converted from quaternion number field to complex number field. An alternating minimization method is applied to solve the model, which is easily implemented and has low computational complexity.
	\item We adopt the rank-decreasing method to estimate the rank of a matrix. Convergence and complexity of the propoesd algorithm is analyzed. Experimental results demonstrate the effectiveness of the propoesd algorithm in color image recovery. Moreover, we compare the propoesd approach with several state-of-the-art tensor-based ones. The results validate the competitive performance of the proposed algorithm compared with the other approaches.
\end{itemize}

The remainder of this paper is organized as follows. Section \ref{sec2} introduces some notations and preliminaries for quaternion algebra.  Section \ref{sec3} reviews the matrix completion theory and proposes our quaternion-based matrix completion model. The detailed overview of quaternion matrix completion algorithm is presented in Section IV. Section V provides some experiments to illustrate the performance of our algorithm,
and compare it with some state-of-the-art methods. Finally, some conclusions are drawn in Section VI.

\section{Notations and preliminaries}
\label{sec2}
In this section, we first summarize some main notations and
then introduce some basic knowledge of quaternion algebra.

\subsection{Notations}
In this paper, $\mathbb{R}$, $\mathbb{C}$ and $\mathbb{H}$ respectively denote the set of real numbers, the set of complex numbers and the set of quaternions. A scalar, a vector, a matrix, and a tensor are written as $a$, $\mathbf{a}$, $\mathbf{A}$, and $\mathcal{A}$, respectively. For
a tensor $\mathcal{A}$, we use the Matlab notation $\mathcal{A}(:,:,k)$  to denote its $k$-th frontal slice. $\ddot{a}$,  $\ddot{\mathbf{a}}$ and $\ddot{\mathbf{A}}$
respectively represent a quaternion scalar, a quaternion vector and a quaternion matrix. $(\cdot)^{\ast}$, $(\cdot)^{-1}$, $(\cdot)^{\dagger}$, $(\cdot)^{T}$ and $(\cdot)^{H}$ denote the
conjugation,  inverse,  Moore-Penrose inverse, transpose and  conjugate transpose respectively. $|\cdot|$, $\|\cdot\|_{1}$, $\|\cdot\|_{F}$ and $\|\cdot\|_{\ast}$ are respectively the absolute value or modulus, the $l_{1}$ norm, the Frobenius norm and the nuclear norm\footnote{For (quaternion) matrix, the nuclear norm is defined as the sum of its singular values. For tensor, the nuclear norm is defined as the nuclear
norms of mode matrices.}. $\circ$ denotes the inner product operation.
${\rm{tr}}\{\cdot\}$ and ${\rm{rank}}(\cdot)$ denote the trace and rank operators respectively. $\mathbf{I}_{K}$ represents the identity matrix of size $K\times K$. And we denote ${\rm{diag}}(\mathbf{x})$ as a diagonal matrix whose diagonal elements are as same as those of $\mathbf{x}$.

\subsection{Basic knowledge of quaternion algebra}
As a natural extension of the complex space, the quaternion space was first introduced by W. Hamilton \cite{articleHamilton84} in 1843. A quaternion $\ddot{q}\in\mathbb{H}$ is composed of a
real component and three imaginary components.
\begin{equation}
\label{equ2}
\ddot{q}=q_{0}+q_{1}i+q_{2}j+q_{3}k,
\end{equation}
where $q_{l}\in\mathbb{R}\: (l=0,1,2,3)$ are real coefficients, $i, j, k$ are
imaginary number units and obey the quaternion rules that
\begin{align}
\left\{
\begin{array}{c}
i^{2}=j^{2}=k^{2}=ijk=-1,\\
ij=-ji=k,\\
jk=-kj=i, \\
ki=-ik=j.
\end{array}
\right.
\end{align}
If the real component $q_{0}=0$, $\ddot{q}$ is named a pure quaternion. Every quaternion $\ddot{q}=q_{0}+q_{1}i+q_{2}j+q_{3}k$ can be uniquely represented as $\ddot{q}=q_{0}+q_{1}i+(q_{2}+q_{3}i)j=c_{1}+c_{2}j$, where $c_{1}$ and $c_{2}$ are complex numbers.

The conjugate and the modulus of a quaternion $\ddot{q}$ are,
respectively, defined as follows:

\begin{align}
\ddot{q}^{\ast}=q_{0}-q_{1}i-q_{2}j-q_{3}k,\\
|\ddot{q}|=\sqrt{q_{0}^{2}+q_{1}^{2}+q_{2}^{2}+q_{3}^{2}}.
\end{align}
Unlike complex number systems, the product of two quaternions $\ddot{q}_{1}$ and $\ddot{q}_{2}$
is noncommunicative, \emph{i.e.}, $\ddot{q}_{1}\ddot{q}_{2}\neq \ddot{q}_{2}\ddot{q}_{1}$  in general.

Analogously, a quaternion matrix $\ddot{\mathbf{Q}}=(\ddot{q}_{mn})\in\mathbb{H}^{M\times N}$ is written
as $\ddot{\mathbf{Q}}=\mathbf{Q}_{0}+\mathbf{Q}_{1}i+\mathbf{Q}_{2}j+\mathbf{Q}_{3}k$, where $\mathbf{Q}_{l}\in\mathbb{R}^{M\times N}\: (l=0,1,2,3)$, $\ddot{\mathbf{Q}}$ is named a pure quaternion matrix when $\mathbf{Q}_{0}=\mathbf{0}$. The Frobenius norm of the quaternion matrix is defined as $\|\ddot{\mathbf{Q}}\|_{F}=\sqrt{\sum_{m=1}^{M}\sum_{n=1}^{N}|\ddot{q}_{mn}|^{2}}=\sqrt{{\rm{tr}}\{(\ddot{\mathbf{Q}})^{H}\ddot{\mathbf{Q}}\}}$.

The most common way to study quaternion matrices is to use their complex representation. Given a quaternion matrix $\ddot{\mathbf{Q}}\in\mathbb{H}^{M\times N}$, it can be uniquely expressed as $\ddot{\mathbf{Q}}=\mathbf{Q}_{a}+\mathbf{Q}_{b}j$, where $\mathbf{Q}_{a},
\mathbf{Q}_{b}\in\mathbb{C}^{M\times N}$. We define the operator $f:\mathbb{H}^{M\times N}\longrightarrow\mathbb{C}^{2M\times 2N}$,
then the complex representation matrix of $\ddot{\mathbf{Q}}=\mathbf{Q}_{a}+\mathbf{Q}_{b}j\in\mathbb{H}^{M\times N}$ is denoted as follows \cite{DBLP:journals/sigpro/BihanM04}:
\begin{equation}
\label{definef}
f(\ddot{\mathbf{Q}})=\left(\begin{array}{cc}
               \mathbf{Q}_{a}	&  \mathbf{Q}_{b}\\
              -\mathbf{Q}^{\ast}_{b}	& \mathbf{Q}^{\ast}_{a}
                    \end{array}\right),
\end{equation}
$f(\ddot{\mathbf{Q}})$ is uniquely determined by $\ddot{\mathbf{Q}}$. Denote $f^{-1}$ as the inverse operator of $f$. Note that if $\mathbf{Q}$ is a complex matrx or real matrx, then
\begin{equation}
\label{definef_c}
f(\mathbf{Q})=\left(\begin{array}{cc}
\mathbf{Q}	&  \mathbf{0}\\
\mathbf{0}	& \mathbf{Q}^{\ast}
\end{array}\right)
\end{equation}
or
\begin{equation}
\label{definef_r}
f(\mathbf{Q})=\left(\begin{array}{cc}
\mathbf{Q}	&  \mathbf{0}\\
\mathbf{0}	& \mathbf{Q}
\end{array}\right),
\end{equation} respectively.

There are some properties of $f$ (\emph{see} Theorem \ref{theorem1} and Theorem \ref{theorem2}).
\begin{theorem}
\label{theorem1}
\cite{10029950538} Let $\ddot{\mathbf{P}}\in\mathbb{H}^{M\times N}$, $\ddot{\mathbf{Q}}\in\mathbb{H}^{M\times N}$, then
\begin{enumerate}
	\item $f(\ddot{\mathbf{P}}\ddot{\mathbf{Q}})=f(\ddot{\mathbf{P}})f(\ddot{\mathbf{Q}})$,
	\item $f(\ddot{\mathbf{P}}+\ddot{\mathbf{Q}})=f(\ddot{\mathbf{P}})+f(\ddot{\mathbf{Q}})$,
	\item $f(\ddot{\mathbf{P}}^{\ast})=f(\ddot{\mathbf{P}})^{\ast}$,
	\item $f(\ddot{\mathbf{P}}^{-1})=f(\ddot{\mathbf{P}})^{-1}$, if $\mathbf{P}^{-1}$ exists,
	\item $\|f(\ddot{\mathbf{P}})\|_{F}^{2}=2\|\ddot{\mathbf{P}}\|_{F}^{2}$,
	\item $\ddot{\mathbf{P}}$  is unitary, Hermitian, or normal if and only if $f(\ddot{\mathbf{P}})$ is unitary, Hermitian, or normal, respectively.
\end{enumerate}
\end{theorem}

\begin{theorem}
\label{theorem2}
Let $\ddot{\mathbf{P}}\in\mathbb{H}^{M\times N}$, we have ${\rm{rank}}(\ddot{\mathbf{P}})=\frac{1}{2}{\rm{rank}}(f(\ddot{\mathbf{P}}))$.
\end{theorem}
The proof of Theorem \ref{theorem2} can be found in Appendix \ref{appendices1}.
Readers can find more details on quaternion algebra in \cite{10029950538, Girard2007Quaternions, Altmann1986Rotations}.

\section{Problem formulation}
\label{sec3}
In this section, we first review the matrix completion theory and then propose our quaternion-based matrix completion model.

\subsection{Matrix completion theory}
Matrix completion problem consists of recovering a matrix from a subset of its entries. The usual structural
assumption on a matrix that makes the problem well posed is that the matrix is low-rank or approximate low-rank. The optimization model for matrix completion was proposed firstly in \cite{Cand2009Exact}, and can be formulated as:
\begin{equation}
\label{equ4}
\begin{split}
&\mathop{{\rm{minimize}}}\limits_{\mathbf{X}}\quad {\rm{rank}}(\mathbf{X})\\
&\text{subject to} \quad  \mathcal{P}_{\Omega}(\mathbf{X}-\mathbf{T})=\mathbf{0},
\end{split}
\end{equation}
where $\mathbf{X}$ is a completed output matrix, $\mathbf{T}$ is an incomplete input
matrix and the $\Omega$ is the entries set, more concretely, if $\mathbf{X}_{mn}$ is observed, then
$(m,n)\in\Omega$, and $\mathcal{P}_{\Omega}$ is the unitary projection onto the linear space of matrices supported on $\Omega$, defined as
\begin{equation*}
(\mathcal{P}_{\Omega}(\mathbf{X}))_{mn}=\left\{
\begin{array}{c}
\!\!\!\mathbf{X}_{mn},\qquad (m,n)\in \Omega, \\
0,\qquad\quad\:  (m,n)\notin\Omega.
\end{array}
\right.
\end{equation*}
Because such rank minimization problem (\ref{equ4}) is generally NP-hard \cite{DBLP:journals/siammax/GillisG11}, various heuristics algorithms have been developed to solve this problem. These methods could be divided into two main categories: nuclear norm minimization method, see \emph{e.g.} \cite{DBLP:journals/mp/MaGC11, DBLP:journals/jmlr/Recht11, DBLP:journals/focm/CandesR09} and low-rank matrix decomposition approach, see \emph{e.g.} \cite{DBLP:conf/stoc/JainNS13, DBLP:journals/mpc/WenYZ12}.

The matrix nuclear norm is a popolar convex surrogate of non-convex rank function, and its minimization method is widely used in practice:
\begin{equation}
\label{equ5}
\begin{split}
&\mathop{{\rm{minimize}}}\limits_{\mathbf{X}}\quad \|\mathbf{X}\|_{\ast}\\
&\text{subject to} \quad  \mathcal{P}_{\Omega}(\mathbf{X}-\mathbf{T})=\mathbf{0}.
\end{split}
\end{equation}
The nuclear norm minimization problem (\ref{equ5}) is
generally solved iteratively in which singular value decompositions (SVD) is involved at each
iteration. So the nuclear norm minimization methods bear the computational cost required by SVD which becomes increasingly expensive as the sizes of the matrices increase \cite{DBLP:journals/mpc/WenYZ12, DBLP:journals/tip/ZhouLLZ18}. Hence, a non-SVD approach, \emph{i.e.} low-rank matrix decomposition, has been proposed in order to more efficiently solve large-scale matrix completion problems.

The low-rank matrix decomposition-based completion problem is formulated in the form of the following optimization problem \cite{DBLP:conf/stoc/JainNS13}:
\begin{equation}
\label{equ6}
\begin{split}
&\mathop{{\rm{minimize}}}\limits_{\mathbf{U}, \mathbf{V}, \mathbf{X}}\quad \frac{1}{2}\|\mathbf{U}\mathbf{V}-\mathbf{X}\|^{2}_{F}\\
&\text{subject to} \quad  \mathcal{P}_{\Omega}(\mathbf{X}-\mathbf{T})=\mathbf{0},
\end{split}
\end{equation}
where $\mathbf{U}\in\mathbb{C}^{M\times K}$, $\mathbf{V}\in\mathbb{C}^{K\times N}$, $\mathbf{X}\in\mathbb{C}^{M\times N}$, and the integer $K$ is the rank
of matrix $\mathbf{X}$.

\subsection{Proposed formulation of quaternion matrix completion}
Quaternion matrix completion can be regard as the generalization of the traditional matrix completion in the quaternion number field, which is to fill in the missing values of a quaternion matrix $\ddot{\mathbf{X}}\in\mathbb{H}^{M\times N}$ under a given subset $\Omega$ of its entries $\{\ddot{\mathbf{X}}_{m,n}|(m,n)\in\Omega \}$. Motivated by traditional matrix completion techniques, low-rank decomposition and nuclear norm minimization, we combine the two approaches and propose our quaternion matrix completion model.
Before that, we first present the following theorem:
\begin{theorem}
	\label{proposition2}
Suppose that $\ddot{\mathbf{X}}\in\mathbb{H}^{M\times N}$, $\ddot{\mathbf{P}}\in\mathbb{H}^{M\times N}$ and $\ddot{\mathbf{Q}}\in\mathbb{H}^{N\times M}$ are three arbitrary quaternion matrices. Then, we have the following properties:
\begin{enumerate}
	\item [(1)] If {\rm{rank}}$(\ddot{\mathbf{X}})=K$, then there extists two quaternion matrices $\ddot{\mathbf{U}}\in\mathbb{H}^{M\times K}$ and $\ddot{\mathbf{V}}\in\mathbb{H}^{K\times N}$ such that
\begin{equation*}
\ddot{\mathbf{X}}=\ddot{\mathbf{U}}\ddot{\mathbf{V}},
\end{equation*}
and they satisfy
\begin{equation*}
{\rm{rank}}(\ddot{\mathbf{U}})={\rm{rank}}(\ddot{\mathbf{V}})=K;
\end{equation*}
 \item [(2)]
${\rm{rank}}(\ddot{\mathbf{P}}\ddot{\mathbf{Q}})\leq {\rm{min}}({\rm{rank}}(\ddot{\mathbf{P}}),  {\rm{rank}}(\ddot{\mathbf{Q}}))$;

\item [(3)] Assume $\ddot{\mathbf{X}}=\ddot{\mathbf{U}}\ddot{\mathbf{V}}$ is a completed output quaternion matrix, $\ddot{\mathbf{T}}$ is an incomplete input
quaternion matrix with rank $K_{0}\leq K$. The nuclear norm minimization problem
\begin{equation}
\label{add_equ1}
\begin{split}
&\mathop{{\rm{minimize}}}\limits_{\ddot{\mathbf{X}}}\quad \|f(\ddot{\mathbf{X}})\|_{\ast}\\
&\text{subject to} \quad  \mathcal{P}_{\Omega}(\ddot{\mathbf{X}}-\ddot{\mathbf{T}})=\mathbf{0}
\end{split}
\end{equation}
is equivalent to the following quadratic optimization problem:
\begin{equation}
\label{add_equ2}
\begin{split}
&\mathop{{\rm{minimize}}}\limits_{f(\ddot{\mathbf{U}}), f(\ddot{\mathbf{V}})}\quad \frac{1}{2}\left( \|f(\ddot{\mathbf{U}})\|_{F}^{2}+\|f(\ddot{\mathbf{V}})\|_{F}^{2}\right) \\
&\text{subject to} \quad  \mathcal{P}_{\Omega}(\ddot{\mathbf{X}}-\ddot{\mathbf{T}})=\mathbf{0}.
\end{split}
\end{equation}
\end{enumerate}
\end{theorem}

The proof of Theorem \ref{proposition2} can be found in Appendix \ref{appendices2}.
Thus, based on the properties (1) and (2) in Theorem \ref{proposition2}, similar to the matrix decomposition method, we can adopt a low-rank quaternion matrix decomposition strategy to deal with the large scale quaternion matrix completion problem more efficiently. Furthermore, we also consider the nuclear norm $\|\ddot{\mathbf{X}}\|_{\ast}$ in our model but replaced
by $\|\ddot{\mathbf{U}}\|_{F}^{2}+\|\ddot{\mathbf{V}}\|_{F}^{2}$ according to property (3) in Theorem \ref{proposition2}.

Accordingly,f/; vdyo,du.qe408888i8bsyn meeddf cv nlj6666jc8iyi the low-rank quaternion matrix completion formulation can be written as follows:

\begin{equation}
\label{equ7}
\begin{split}
&\mathop{{\rm{minimize}}}\limits_{\ddot{\mathbf{U}}, \ddot{\mathbf{V}}, \ddot{\mathbf{X}}}\quad \frac{1}{2}\|\ddot{\mathbf{U}}\ddot{\mathbf{V}}-\ddot{\mathbf{X}}\|^{2}_{F}+\frac{\lambda}{2}\left( \|\ddot{\mathbf{U}}\|_{F}^{2}+\|\ddot{\mathbf{V}}\|_{F}^{2}\right) \\
&\text{subject to} \quad  \mathcal{P}_{\Omega}(\ddot{\mathbf{X}}-\ddot{\mathbf{T}})=\mathbf{0},
\end{split}
\end{equation}
where $\lambda$ is a nonnegative parameter.

\section{Proposed algorithm}
In this section, we first show how to solve the optimization
problem (\ref{equ7}), then we introduce a rank-decreasing method to adjust the rank of a matrix. Finally, we provide the convergence and complexity analyses of the proposed algorithm.

\subsection{Optimization process}
On account of the noncommutativity of the multiplication in quaternion space,
the definition and computation of the gradient of quaternion matrix function are generally much more complicated than those in complex space \cite{DBLP:journals/tsp/XuM15}, which hugely increases
the difficulty to handle the quaternion-based optimization problems. Therefore, based on the defined operator $f$ in (\ref{definef}) and its
properties in Theorem \ref{theorem1}, we tend to convert the problem (\ref{equ7}) to that in the complex number field and reformulate as follows, which is differentiable and separable among
its blocks:
\begin{equation}
\label{equ8}
\begin{split}
&\mathop{{\rm{minimize}}}\limits_{f(\ddot{\mathbf{U}}), f(\ddot{\mathbf{V}}), \ddot{\mathbf{X}}} \quad \frac{1}{2}\|f(\ddot{\mathbf{U}})f(\ddot{\mathbf{V}})-f(\ddot{\mathbf{X}})\|^{2}_{F} \\
&\qquad\qquad\qquad\quad +\frac{\lambda}{2}\left( \|f(\ddot{\mathbf{U}})\|_{F}^{2}+\|f(\ddot{\mathbf{V}})\|_{F}^{2}\right)\\
&\text{subject to} \quad  \mathcal{P}_{\Omega}(\ddot{\mathbf{X}}-\ddot{\mathbf{T}})=\mathbf{0}.
\end{split}
\end{equation}
Note that $f(\ddot{\mathbf{U}})\in\mathbb{C}^{2M\times 2K}$ and $f(\ddot{\mathbf{V}})\in\mathbb{C}^{2K\times 2N}$ are all complex-valued matrices.

Although, it is obvious that problem (\ref{equ8}) is non-convex itself,  it
is convex with respect to each single variable. Hence, we adopt a simple but efficient
iterative scheme to solve the optimization problem (\ref{equ8}) by using an alternating minimization approach. More specifically, we update only one of the variables $f(\ddot{\mathbf{U}})$, $f(\ddot{\mathbf{V}})$ and $\ddot{\mathbf{X}}$ each time while remaining the other two fixed, and three variables all will be updated sequentially and iteratively.

Letting
\begin{align}
\label{ob_function}
\mathcal{G}(f(\ddot{\mathbf{U}}), f(\ddot{\mathbf{V}}), \ddot{\mathbf{X}})=&\frac{1}{2}\|f(\ddot{\mathbf{U}})f(\ddot{\mathbf{V}})-f(\ddot{\mathbf{X}})\|^{2}_{F} \nonumber\\
& +\frac{\lambda}{2}\left( \|f(\ddot{\mathbf{U}})\|_{F}^{2}+\|f(\ddot{\mathbf{V}})\|_{F}^{2}\right),
\end{align}
we perform the updates as
\begin{subequations}
\label{equ9}
\begin{align}
&f(\ddot{\mathbf{U}})^{\tau+1}=\mathop{{\rm{arg\, min}}}\limits_{f(\ddot{\mathbf{U}})}\:\mathcal{G}(f(\ddot{\mathbf{U}}), f(\ddot{\mathbf{V}})^{\tau}, \ddot{\mathbf{X}}^{\tau}),\\
&f(\ddot{\mathbf{V}})^{\tau+1}=\mathop{{\rm{arg\, min}}}\limits_{f(\ddot{\mathbf{V}})}\:\mathcal{G}(f(\ddot{\mathbf{U}})^{\tau+1}, f(\ddot{\mathbf{V}}), \ddot{\mathbf{X}}^{\tau}),\\
&\ddot{\mathbf{X}}^{\tau+1}=\mathop{{\rm{arg\, min}}}\limits_{ \mathcal{P}_{\Omega}(\ddot{\mathbf{X}}-\ddot{\mathbf{T}})=\mathbf{0}}\:\mathcal{G}(f(\ddot{\mathbf{U}})^{\tau+1}, f(\ddot{\mathbf{V}})^{\tau+1}, \ddot{\mathbf{X}}),
\end{align}
\end{subequations}
where $\tau$ is the iteration index.

By introducing a Lagrange multiplier $\mathbf{\Upsilon}$ for the constraint $\mathcal{P}_{\Omega}(\ddot{\mathbf{X}}-\ddot{\mathbf{T}})=\mathbf{0}$, the
Lagrangian function of (\ref{equ8}) is defined as
\begin{align*}
\mathcal{Q}(f(\ddot{\mathbf{U}}), f(\ddot{\mathbf{V}}), \ddot{\mathbf{X}}, \mathbf{\Upsilon})=&\mathcal{G}(f(\ddot{\mathbf{U}}), f(\ddot{\mathbf{V}}), \ddot{\mathbf{X}})\\
&-\mathbf{\Upsilon}\circ\mathcal{P}_{\Omega}(\ddot{\mathbf{X}}-\ddot{\mathbf{T}}).
\end{align*}
Differentiating the function $\mathcal{Q}(f(\ddot{\mathbf{U}}), f(\ddot{\mathbf{V}}), \ddot{\mathbf{X}}, \mathbf{\Upsilon})$, we have the following Karush-Kuhn-Tucker (KKT) Conditions:
\begin{subequations}
\begin{align}
(f(\ddot{\mathbf{U}})f(\ddot{\mathbf{V}})-f(\ddot{\mathbf{X}}))f(\ddot{\mathbf{V}})^{H}+\lambda f(\ddot{\mathbf{U}}) &=\mathbf{0},\label{equ10_1}\\
f(\ddot{\mathbf{U}})^{H}(f(\ddot{\mathbf{U}})f(\ddot{\mathbf{V}})-f(\ddot{\mathbf{X}}))+\lambda f(\ddot{\mathbf{V}}) &=\mathbf{0},\label{equ10_2}\\
\mathcal{P}_{\Omega^{c}}\left(\ddot{\mathbf{X}}-f^{-1}( f(\ddot{\mathbf{U}})f(\ddot{\mathbf{V}}))\right) &=\mathbf{0},\label{equ10_3}\\
\mathcal{P}_{\Omega}(\ddot{\mathbf{X}}-\ddot{\mathbf{T}})&=\mathbf{0},\label{equ10_4}\\
\mathcal{P}_{\Omega}\left(\ddot{\mathbf{X}}-f^{-1}( f(\ddot{\mathbf{U}})f(\ddot{\mathbf{V}}))\right)-\mathbf{\Upsilon} &=\mathbf{0}.\label{equ10_5}
\end{align}
\end{subequations}
Thus, the updates in (\ref{equ9}) can be explicitly written as follows:

\begin{align}
&f(\ddot{\mathbf{U}})^{\tau+1}=f(\ddot{\mathbf{X}})^{\tau}(f(\ddot{\mathbf{V}})^{\tau})^{H}\Psi_{\ddot{\mathbf{V}}},\label{equ11}\\
&f(\ddot{\mathbf{V}})^{\tau+1}=\Phi_{\ddot{\mathbf{U}}}(f(\ddot{\mathbf{U}})^{\tau+1})^{H}f(\ddot{\mathbf{X}})^{\tau},\label{equ12}
\end{align}
where
\begin{align*}
&\Psi_{\ddot{\mathbf{V}}}=\left( f(\ddot{\mathbf{V}})^{\tau}(f(\ddot{\mathbf{V}})^{\tau})^{H}+\lambda\mathbf{I}_{2K}\right)^{\dagger},\\
&\Phi_{\ddot{\mathbf{U}}}=\left( (f(\ddot{\mathbf{U}})^{\tau+1})^{H}f(\ddot{\mathbf{U}})^{\tau+1}+\lambda\mathbf{I}_{2K}\right)^{\dagger}.
\end{align*}
Then, we can directly obtain $\ddot{\mathbf{X}}^{\tau+1}$ as
\begin{equation}
\label{equ14}
\ddot{\mathbf{X}}^{\tau+1}=\mathcal{P}_{\Omega^{c}}\left( f^{-1}( f(\ddot{\mathbf{U}}^{\tau+1})f(\ddot{\mathbf{V}}^{\tau+1}))\right) + \ddot{\mathbf{T}},
\end{equation}
where $\Omega^{c}$ is the complement of $\Omega$, and we have used the fact that $\mathcal{P}_{\Omega^{c}}(\ddot{\mathbf{T}})=\mathbf{0}$ in (\ref{equ14}).

\subsection{Rank estimation based on a rank-decreasing method}
A proper estimation to the rank $2K$ (labeled as r) for the model (\ref{equ8}) is essential for the success of the proposed algorithm. Although the target rank can be adjusted manually, it would be usually
time-consuming for large-scale data. We introduce a rank-decreasing method to estimate the rank, which is similar to that in \cite{DBLP:journals/mpc/WenYZ12}.

This method starts from an input overestimated rank r of $f(\ddot{\mathbf{X}})$, \emph{i.e.}, ${\rm{r}}>{\rm{rank}}(f(\ddot{\mathbf{X}}))$.
Suppose that the rank of $f(\ddot{\mathbf{X}})^{\tau}$ is ${\rm{r}}^{\tau}$. We compute the eigenvalues of $( f(\ddot{\mathbf{U}})^{\tau})^{H}f(\ddot{\mathbf{U}})^{\tau}\mathbf{\Pi}^{\tau}$, where $\mathbf{\Pi}^{\tau}$ is a permutation matrix so that all
these eigenvalues are ordered non-increasing, \emph{i.e.}, $\mathbf{d}^{\tau}_{1}\geq\mathbf{d}^{\tau}_{2}\geq, \ldots, \geq \mathbf{d}^{\tau}_{{\rm{r}}^{\tau}}$. Then, we compute the quotient sequnce $\hat{\mathbf{d}}^{\tau}_{m}=\mathbf{d}^{\tau}_{m}/\mathbf{d}^{\tau}_{m+1},\:(m=1, \ldots, {\rm{r}}^{\tau}-1)$. Assume that
\begin{equation*}
p^{\tau}=\mathop{{\rm{arg\, max}}}\limits_{1\leq m\leq {\rm{r}}^{\tau}-1}\hat{\mathbf{d}}^{\tau}_{m},
\end{equation*}
and define
\begin{equation}
\label{geprank}
\mu^{\tau}=\frac{({\rm{r}}^{\tau}-1)\hat{\mathbf{d}}^{\tau}_{p^{\tau}}}{\sum_{m\neq p^{\tau}}\hat{\mathbf{d}}^{\tau}_{m}}.
\end{equation}
If $\mu^{\tau}\geq 10$, \emph{i.e.}, there being a large drop in the estimated rank of the $f(\ddot{\mathbf{X}})^{\tau}$, we should reduce ${\rm{r}}^{\tau}$ to $p^{\tau}$. Then, assuming
$\mathbf{L}^{\tau}\mathbf{\Sigma}^{\tau}(\mathbf{R}^{H})^{\tau}$  is the SVD of $f(\ddot{\mathbf{U}})^{\tau}f(\ddot{\mathbf{V}})^{\tau}$, we can update $f(\ddot{\mathbf{U}})^{\tau}=\mathbf{L}^{\tau}_{p^{\tau}}\mathbf{\Sigma}^{\tau}_{p^{\tau}}$ and $f(\ddot{\mathbf{V}})^{\tau}=(\mathbf{R}^{H})^{\tau}_{p^{\tau}}$, where $\mathbf{L}^{\tau}_{p^{\tau}}$ consists of the first $p^{\tau}$ columns of $f(\ddot{\mathbf{U}})^{\tau}$, and $(\mathbf{R}^{H})^{\tau}_{p^{\tau}}$ and $\mathbf{\Sigma}^{\tau}_{p^{\tau}}$ are obtained accordingly. Note that, doing only one time of this rank-adjusting scheme is generally enough during the whole iterative process. Hence, the computational complexity generated from the SVD of $f(\ddot{\mathbf{U}})^{\tau}f(\ddot{\mathbf{V}})^{\tau}$ in this rank-adjusting process is negligible relative to that of the whole iterative process.

Finally, the proposed \textbf{L}ow-\textbf{R}ank \textbf{Q}uaternion \textbf{M}atrix \textbf{C}ompletion (LRQMC) algorithm can be summarized as shown in TABLE \ref{tab_algorithm}.

\begin{table}[htbp]
	\caption{The low-rank quaternion matrix completion (LRQMC) algorithm.}
	\hrule
	\label{tab_algorithm}
	\begin{algorithmic}[1]
		\REQUIRE   The quaternion matrix data $\ddot{\mathbf{X}}\in\mathbb{H}^{M\times N}$, the observed set $\Omega$, and the initialized rank ${\rm{r}}^{0}$.
		\STATE \textbf{Initialize} $f(\ddot{\mathbf{U}})^{0}\in\mathbb{C}^{2M\times {\rm{r}}^{0}}$ and $f(\ddot{\mathbf{V}})^{0}\in\mathbb{C}^{{\rm{r}}^{0}\times 2N}$ randomly, and appropriate parameter $\lambda>0$.
		\STATE \textbf{Repeat}
		\STATE Fix $f(\ddot{\mathbf{V}})^{\tau}$ and $f(\ddot{\mathbf{X}})^{\tau}$ to update $f(\ddot{\mathbf{U}})^{\tau+1}$ by (\ref{equ11}), \emph{i.e.}, $f(\ddot{\mathbf{U}})^{\tau+1}\longleftarrow f(\ddot{\mathbf{X}})^{\tau}(f(\ddot{\mathbf{V}})^{\tau})^{H}\Psi_{\ddot{\mathbf{V}}}$.
		\STATE Fix $f(\ddot{\mathbf{U}})^{\tau+1}$ and $f(\ddot{\mathbf{X}})^{\tau}$ to update $f(\ddot{\mathbf{V}})^{\tau+1}$ by (\ref{equ12}), \emph{i.e.},
		$f(\ddot{\mathbf{V}})^{\tau+1}\longleftarrow\Phi_{\ddot{\mathbf{U}}}(f(\ddot{\mathbf{U}})^{\tau+1})^{H}f(\ddot{\mathbf{X}})^{\tau}$.
		\STATE Fix $f(\ddot{\mathbf{U}})^{\tau+1}$ and $f(\ddot{\mathbf{V}})^{\tau+1}$ to update $\ddot{\mathbf{X}}^{\tau+1}$ by (\ref{equ14}), \emph{i.e.},
		$\ddot{\mathbf{X}}^{\tau+1}\longleftarrow \mathcal{P}_{\Omega^{c}}\left( f^{-1}( f(\ddot{\mathbf{U}}^{\tau+1})f(\ddot{\mathbf{V}}^{\tau+1}))\right) + \ddot{\mathbf{T}}$.
		\IF{$\mu^{\tau}\geq 10$ in (\ref{geprank})}
		\STATE Apply rank-decreasing method to adjust ${\rm{r}}^{\tau}$ and the sizes of $f(\ddot{\mathbf{U}})^{\tau+1}$ and $f(\ddot{\mathbf{V}})^{\tau+1}$.
		\ENDIF
		\STATE  $\tau\longleftarrow \tau+1$.
		\STATE \textbf{Until convergence}
		\ENSURE   $f(\ddot{\mathbf{U}})^{\tau+1}$, $f(\ddot{\mathbf{V}})^{\tau+1}$ and $\ddot{\mathbf{X}}^{\tau+1}$
	\end{algorithmic}
	\hrule
\end{table}

\subsection{The convergence and the computational complexity analyses}

\textbf{Convergence analysis:}
From (\ref{equ10_5}), we can clearly see that $\mathcal{P}_{\Omega}\left(\ddot{\mathbf{X}}-f^{-1}( f(\ddot{\mathbf{U}})f(\ddot{\mathbf{V}}))\right)=\mathbf{\Upsilon}$, \emph{i.e.}, the multiplier matrix $\mathbf{\Upsilon}$ measures the residual $\ddot{\mathbf{X}}-f^{-1}( f(\ddot{\mathbf{U}})f(\ddot{\mathbf{V}}))$ in $\Omega$ and thus has no effect in the process of determining $f(\ddot{\mathbf{U}})$, $f(\ddot{\mathbf{V}})$ and $\ddot{\mathbf{X}}$. Therefore, for simplicity, we just discuss $\mathcal{G}(f(\ddot{\mathbf{U}}), f(\ddot{\mathbf{V}}), \ddot{\mathbf{X}})$ in (\ref{ob_function}).

Since the Hessian matrices of $\mathcal{G}(f(\ddot{\mathbf{U}}), f(\ddot{\mathbf{V}}), \ddot{\mathbf{X}})$ \emph{ w.r.t.} $f(\ddot{\mathbf{U}})$ and $ f(\ddot{\mathbf{V}})$ are respectively $f(\ddot{\mathbf{V}})f(\ddot{\mathbf{V}})^{H}+\lambda\mathbf{I}_{2K}$ and $f(\ddot{\mathbf{U}})f(\ddot{\mathbf{U}})^{H}+\lambda\mathbf{I}_{2K}$ which are positive semidefinite matrices (they are even positive definite when $\lambda>\mathbf{0}$). Hence, for any $\tau\geq 0$, we have $\mathcal{G}(f(\ddot{\mathbf{U}})^{\tau}, f(\ddot{\mathbf{V}})^{\tau}, \ddot{\mathbf{X}}^{\tau})-\mathcal{G}(f(\ddot{\mathbf{U}})^{\tau+1}, f(\ddot{\mathbf{V}})^{\tau+1}, \ddot{\mathbf{X}}^{\tau})\geq 0$. On the other hand, we note that $\ddot{\mathbf{X}}^{\tau+1}$ is the optimal solution to problem (\ref{equ7}):
\begin{align}
\label{convergence_1}
\ddot{\mathbf{X}}^{\tau+1}=&\mathop{{\rm{arg\, min}}}\limits_{ \mathcal{P}_{\Omega}(\ddot{\mathbf{X}}-\ddot{\mathbf{T}})=\mathbf{0}}\frac{1}{2}\|\ddot{\mathbf{U}}^{\tau+1}\ddot{\mathbf{V}}^{\tau+1}-\ddot{\mathbf{X}}\|^{2}_{F}\nonumber\\
&+\frac{\lambda}{2}\left( \|\ddot{\mathbf{U}}^{\tau+1}\|_{F}^{2}+\|\ddot{\mathbf{V}}^{\tau+1}\|_{F}^{2}\right).
\end{align}
At the same time, we note that $\mathcal{P}_{\Omega}(\ddot{\mathbf{X}}^{\tau+1}-\ddot{\mathbf{T}})=\mathbf{0}$, \emph{i.e.}, $\ddot{\mathbf{X}}^{\tau+1}$ is a feasible solution to problem (\ref{convergence_1}). So the inequality
$\frac{1}{2}\|\ddot{\mathbf{U}}^{\tau+1}\ddot{\mathbf{V}}^{\tau+1}-\ddot{\mathbf{X}}^{\tau+1}\|^{2}_{F}\leq \frac{1}{2}\|\ddot{\mathbf{U}}^{\tau+1}\ddot{\mathbf{V}}^{\tau+1}-\ddot{\mathbf{X}}^{\tau}\|^{2}_{F}$ holds,
\emph{i.e.}, the inequality
$\frac{1}{2}\|f(\ddot{\mathbf{U}})^{\tau+1}f(\ddot{\mathbf{V}})^{\tau+1}-f(\ddot{\mathbf{X}})^{\tau+1}\|^{2}_{F}\leq \frac{1}{2}\|f(\ddot{\mathbf{U}})^{\tau+1}f(\ddot{\mathbf{V}})^{\tau+1}-f(\ddot{\mathbf{X}})^{\tau}\|^{2}_{F}$ holds. Hence, for any $\tau\geq 0$, we can obtain $\mathcal{G}(f(\ddot{\mathbf{U}})^{\tau+1}, f(\ddot{\mathbf{V}})^{\tau+1}, \ddot{\mathbf{X}}^{\tau})-\mathcal{G}(f(\ddot{\mathbf{U}})^{\tau+1}, f(\ddot{\mathbf{V}})^{\tau+1}, \ddot{\mathbf{X}}^{\tau+1})\geq 0$. Then, it follows that
\begin{align*}
&\mathcal{G}(f(\ddot{\mathbf{U}})^{\tau}, f(\ddot{\mathbf{V}})^{\tau}, \ddot{\mathbf{X}}^{\tau})-\mathcal{G}(f(\ddot{\mathbf{U}})^{\tau+1}, f(\ddot{\mathbf{V}})^{\tau+1}, \ddot{\mathbf{X}}^{\tau+1})\\
&=\mathcal{G}(f(\ddot{\mathbf{U}})^{\tau}, f(\ddot{\mathbf{V}})^{\tau}, \ddot{\mathbf{X}}^{\tau})-\mathcal{G}(f(\ddot{\mathbf{U}})^{\tau+1}, f(\ddot{\mathbf{V}})^{\tau+1}, \ddot{\mathbf{X}}^{\tau})\\
&\quad+\mathcal{G}(f(\ddot{\mathbf{U}})^{\tau+1}, f(\ddot{\mathbf{V}})^{\tau+1}, \ddot{\mathbf{X}}^{\tau})\\
&\quad-\mathcal{G}(f(\ddot{\mathbf{U}})^{\tau+1}, f(\ddot{\mathbf{V}})^{\tau+1}, \ddot{\mathbf{X}}^{\tau+1})\\
&\geq 0.
\end{align*}
Consequently, the function $\mathcal{G}(f(\ddot{\mathbf{U}})^{\tau}, f(\ddot{\mathbf{V}})^{\tau}, \ddot{\mathbf{X}}^{\tau})$ decreases monotonically, and it is obvious that $\mathcal{G}(f(\ddot{\mathbf{U}})^{\tau}, f(\ddot{\mathbf{V}})^{\tau}, \ddot{\mathbf{X}}^{\tau})\geq 0$, so the theoretical convergence of the proposed algorithm (LRQMC) can
be guaranteed.

\textbf{Computational complexity:} We analyze the computational complexity within one iteration for LRQMC algorithm provided in TABLE \ref{tab_algorithm}. When
updating $f(\ddot{\mathbf{U}})$ and  $f(\ddot{\mathbf{V}})$ respectively by (\ref{equ11}) and (\ref{equ12}), the computational cost is about $\mathcal{O}\left({\rm{\hat{r}}}({\rm{\hat{r}}}^{2}+N{\rm{\hat{r}}}+M{\rm{\hat{r}}}+MN)\right)$, where ${\rm{\hat{r}}}$ is the estimated rank of $f(\ddot{\mathbf{X}})$. Then the computational cost of updating  $\ddot{\mathbf{X}}$ by (\ref{equ14}) is about $\mathcal{O}(MN{\rm{\hat{r}}})$. In rank estimation process, the estimated rank is detected based on the economy-size QR decomposition whose computational cost is about $\mathcal{O}({\rm{\hat{r}}}^{3})$. Hence, the total computational cost of LRQMC algorithm at each iteration is about $\mathcal{O}\left({\rm{\hat{r}}}({\rm{\hat{r}}}^{2}+N{\rm{\hat{r}}}+M{\rm{\hat{r}}}+MN)\right)$.

\section{Simulation results}

In this section, simulations on some natural color images are conducted to evaluate the performance of the proposed LRQMC algorithm. And we compare it
with several existing state-of-the-art tensor-based methods, including TCTF \cite{DBLP:journals/tip/ZhouLLZ18}, SPC \cite{DBLP:journals/tsp/YokotaZC16}, TMac (including TMac-inc and TMac-dec) \cite{DBLP:journals/corr/XuHYS13}, STDC \cite{DBLP:journals/pami/ChenHL14} and LRTC (including FaLRTC and SiLRTC) \cite{DBLP:journals/pami/LiuMWY13}. All the simulations are run in MATLAB $2014b$ under Windows $10$ on a personal computer with $2.20$GHz CPU and $8$GB memory.

A color image is a 3-way tensor defined by two indices for spatial variables and one index for color mode \cite{DBLP:journals/tip/BenguaPTD17}.  All the images, in our simulation, are initially
represented by 3-way tensors $\mathcal{T}\in \mathbb{R}^{M\times N\times s}$, where $M\times N$ is the number of pixels
in the image and $s=3$ is the number of colors (red, green and blue). For LRQMC algorithm,  each image  is reshaped as a  pure
quaternion matrix $\ddot{\mathbf{T}}\in\mathbb{H}^{M\times N}$ by using the following way:
\begin{equation*}
\ddot{\mathbf{T}}= \mathcal{T}(:,:,1)i+\mathcal{T}(:,:,2)j+\mathcal{T}(:,:,3)k.
\end{equation*}
In addition, we uniformly generate the index set $\Omega$ at Gaussian random distribution, and define the sampling ratio (SR) as:
\begin{equation*}
{\rm{SR}}=\frac{{\rm{numel}}(\Omega)}{M\times N\times s},
\end{equation*}
where ${\rm{numel}}(\Omega)$  represents the number of observation elements in the index set $\Omega$.

\textbf{Quantitative assessment:} In order to evaluate
the performance of proposed algorithm, except visual quality, we employ four quantitative quality indexes, including the relative square error (RSE), the  peak signal-to-noise ratio (PSNR), the structure similarity (SSIM) and the feature similarity (FSIM), which are respectively defined as follows:
\begin{equation*}
{\rm{RSE}}=10{\rm{log}}10\left( \frac{\|\mathcal{X}-\mathcal{T}\|_{F}}{\|\mathcal{T}\|_{F}}\right),
\end{equation*}
where $\mathcal{X}$ and $\mathcal{T}$ are the recovered and truth data, respectively.
\begin{equation*}
{\rm{PSNR}}=10{\rm{log}}10\left( \frac{{\rm{Peakval}}^{2}}{{\rm{MSE}}}\right),
\end{equation*}
where ${\rm{Peakval}}$ is taken from the range of the image datatype (\emph{e.g.}, for uint8 image it is 255), ${\rm{MSE}}$ is the mean square error, \emph{i.e.} ${\rm{MSE}}=\|\mathcal{X}-\mathcal{T}\|_{F}^{2}/{\rm{numel}}(\mathcal{X})$.
\begin{equation*}
{\rm{SSIM}}=\frac{(2\mu_{\mathcal{T}}\mu_{\mathcal{X}}+C_{1})(2\sigma_{\mathcal{T}\mathcal{X}}+C_{2})}{(\mu_{\mathcal{T}}^{2}+\mu_{\mathcal{X}}^{2}+C_{1})(\sigma_{\mathcal{T}}^{2}+\sigma_{\mathcal{X}}^{2}+C_{2})},
\end{equation*}
where $\mu_{\mathcal{T}}$, $\mu_{\mathcal{X}}$, $\sigma_{\mathcal{T}}$, $\sigma_{\mathcal{X}}$ and $\sigma_{\mathcal{T}\mathcal{X}}$ are the local means, standard deviations, and cross-covariance for images $\mathcal{T}$ and $\mathcal{X}$, $C_{1}=(0.01L)^{2}$,  $C_{2}=(0.03L)^{2}$,  $C_{3}=C_{2}/2$, $L$ is the specified dynamic range of the pixel values.
\begin{equation*}
{\rm{FSIM}}=\frac{\sum_{z\in \Delta}S_{L}(z)PC_{m}(z)}{\sum_{z\in \Delta}PC_{m}(z)},
\end{equation*}
where $\Delta$ demotes the whole image spatial domain. The phase
congruency for position $z$ of image $\mathcal{T}$  is denoted as $PC_{x}(\mathcal{T})$, then $PC_{m}(z)={\rm{max}}\{PC_{\mathcal{T}(z)},PC_{\mathcal{X}(z)}\}$, $S_{L}(z)$ isthegradient
magnitude for position $z$.

For LRQMC algorithm, $\mathcal{X}(:,:,1)={\rm{Imag}}_{1}(\ddot{\mathbf{T}})$, $\mathcal{X}(:,:,2)={\rm{Imag}}_{2}(\ddot{\mathbf{T}})$ and $\mathcal{X}(:,:,3)={\rm{Imag}}_{3}(\ddot{\mathbf{T}})$, where ${\rm{Imag}}_{n}(\ddot{\mathbf{T}})\:(n=1,2,3)$ denotes $n$-th image part of $\ddot{\mathbf{T}}$.

\textbf{Datasets:} In the simulations, we use two color image datasets: Berkeley segmentation dataset and Kodak PhotoCD dataset. The statistics of these two datasets are briefly summarized below:
\begin{itemize}
	\item[-]  Berkeley Segmentation Dataset (BSD):\footnote{\url{ https://www2.eecs.berkeley.edu/Research/Projects/CS/vision/bsds/}} There are $300$
	clean color images of size $481\times 321\times 3$ in the whole dataset.
	\item[-] Kodak PhotoCD Dataset (Kodak):\footnote{\url{http://r0k.us/graphics/kodak/}} The whole dataset consists of $24$ clean color images of size $512 \times 768\times 3$. 	
\end{itemize}

\begin{figure*}[htbp]
	\centering
	\subfigure[]{\includegraphics[width=2.5cm,height=2.8cm]{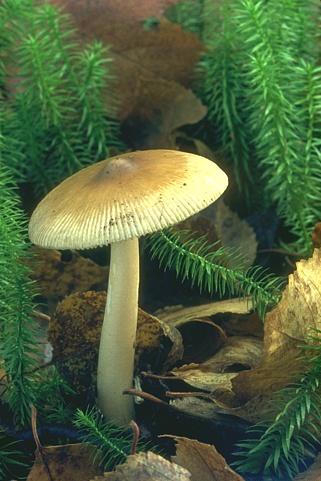}}
	\subfigure[]{\includegraphics[width=2.5cm,height=2.8cm]{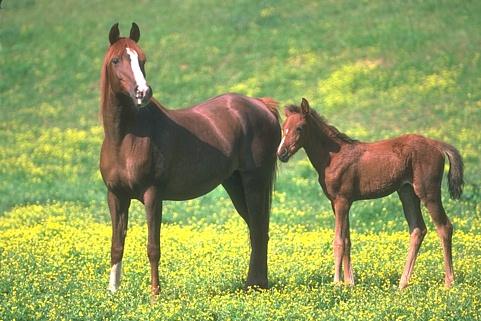}}
	\subfigure[]{\includegraphics[width=4.5cm,height=3.1cm]{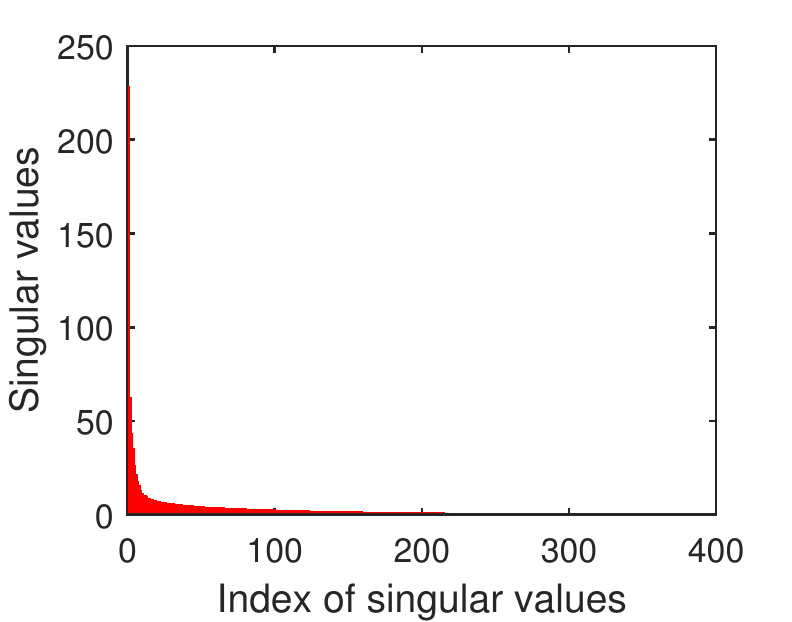}}
	\subfigure[]{\includegraphics[width=4.5cm,height=3.1cm]{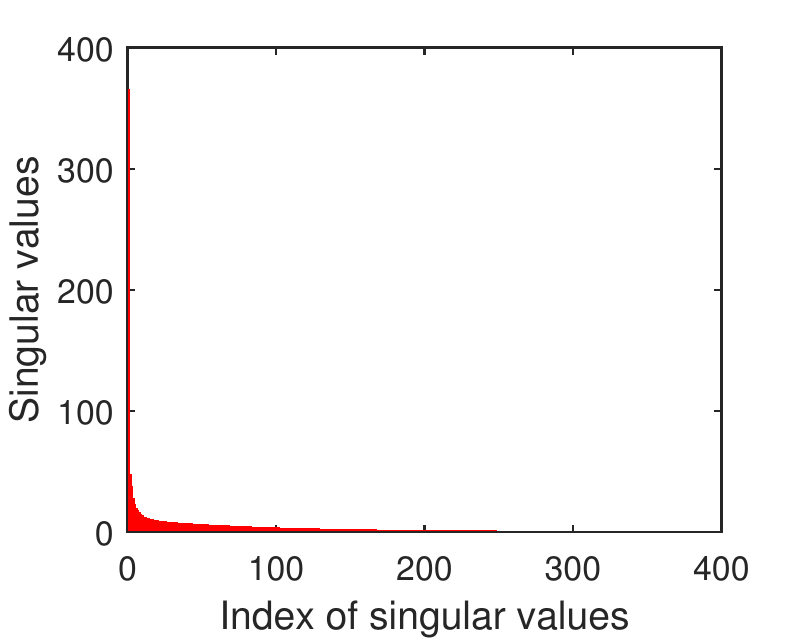}}\\ 
	\subfigure[]{\includegraphics[width=2.5cm,height=2.8cm]{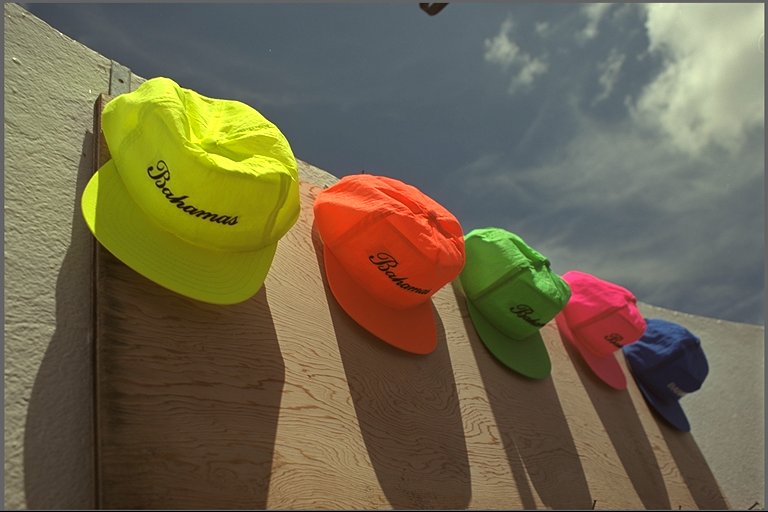}}
	\subfigure[]{\includegraphics[width=2.5cm,height=2.8cm]{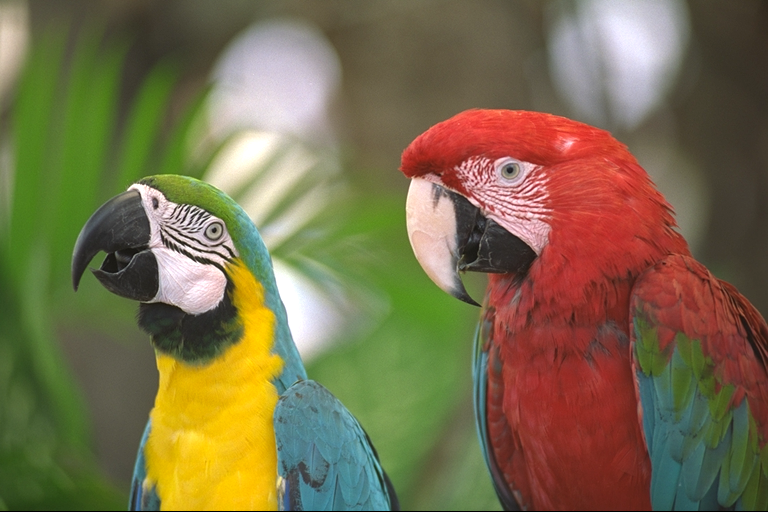}}
	\subfigure[]{\includegraphics[width=4.5cm,height=3.1cm]{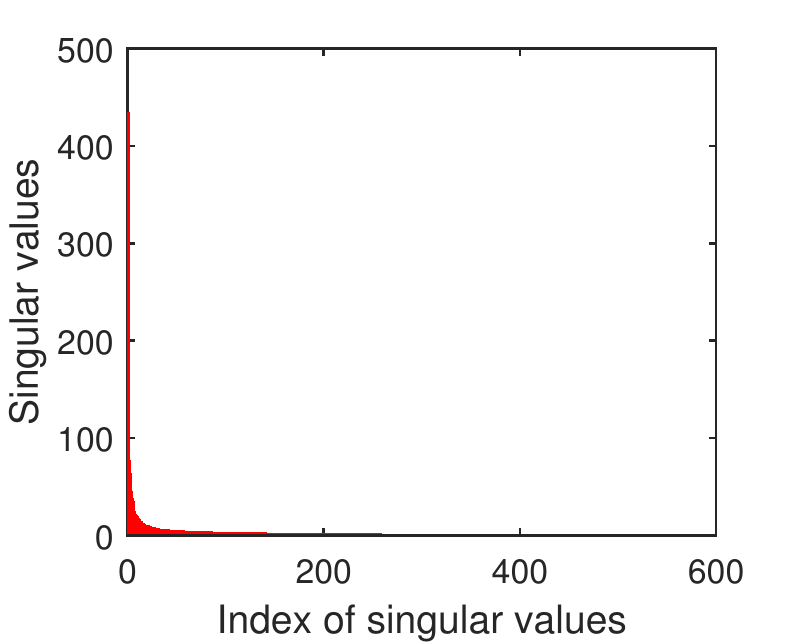}}
	\subfigure[]{\includegraphics[width=4.5cm,height=3.1cm]{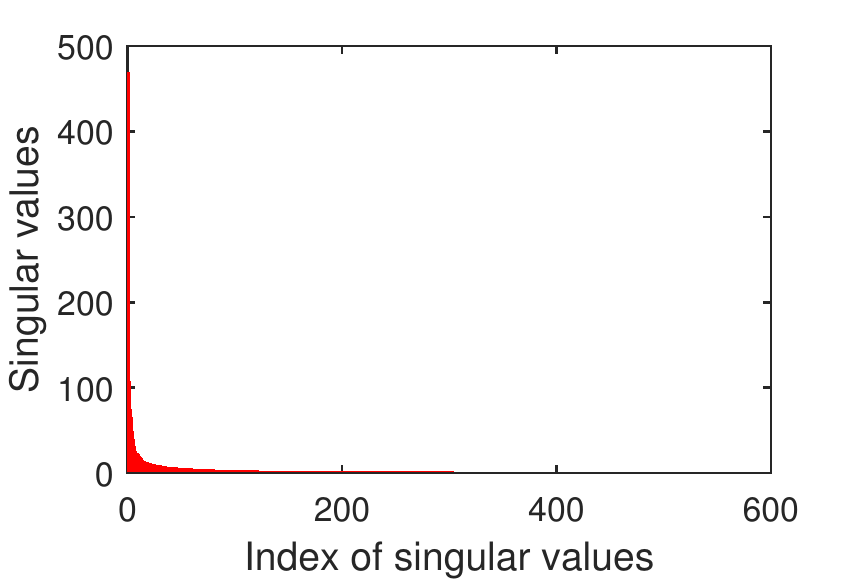}}\\
	
	\caption{Illustration of the low-rank property of the images in the two color image datasets. (a) and (b) are two images randomly selected from the BSD,  (c) and (d) respectively display the singular values of (a) and (b). (e) and (f) are two images randomly selected from the Kodak, (g) and (h) respectively display the singular values of (e) and (f).}
	\label{fig_rank}
\end{figure*}

We first show that these color images can be well approximated by the low-rank quaternion matrices. Actually, as mentioned in \cite{4797640Candes2019,DBLP:journals/tip/ZhouLLZ18}, when the image data is arranged into matrices or tensors, they lie on a
union of low-rank subspaces approximately, which indicate the low-rank structure of the image data. This is also true for quaternion matrices data. For instance, in Fig. \ref{fig_rank} we display the singular values of four images (reshaped as pure quaternion matrices) selected from the two color image datasets randomly. One can obviously see that most of the singular values are very close to $0$, and much smaller than the first several larger singular values. So we could say that these color images can be well approximated by the low-rank quaternion matrices as we desired.

\textbf{Parameter settings:} For our LRQMC algorithm, the initial rank ${\rm{r}}$ of $f(\ddot{\mathbf{X}})$ is set as ${\rm{r}}=50$, and we set  $\lambda=0.5$. For TCTF,
we set the initialized rank ${\rm{r}}^{0}=[30, 30, 30]$ the same as that in \cite{DBLP:journals/tip/ZhouLLZ18}. For SPC, we use QV constraint with $\rho=[1.0, 1.0, 0]$. For TMac-inc, we set the initialized rank ${\rm{r}}^{0}=[3, 3, 3]$ with increment $2$. For TMac-dec, we set the initialized rank ${\rm{r}}^{0}=[30, 30, 30]$. And we set $\alpha_{n}=\frac{1}{3},\: n=1,2,3$ for both TMac-inc and TMac-dec as suggested in \cite{DBLP:journals/corr/XuHYS13}. For SiLRTC, according to \cite{DBLP:journals/pami/LiuMWY13}, the weight parameter $\alpha=\theta/\|\theta\|_{1}$, where $\theta=[1,1,1e^{-3}]$. In addition, the stopping criteria for all the algorithms that we adopted are the difference between the values of $\varepsilon:=\|\mathcal{X}-\mathcal{T}\|_{F}$ in two consecutive iterations, \emph{i.e.} $|\varepsilon^{\tau}-\varepsilon^{\tau+1}|<1e^{-3}$,
where $\tau$ is the iteration index, and the maximum number of iterations is $1,000$.

\begin{figure*}[htbp]
	\centering
		\includegraphics[width=15.8cm,height=12.8cm]{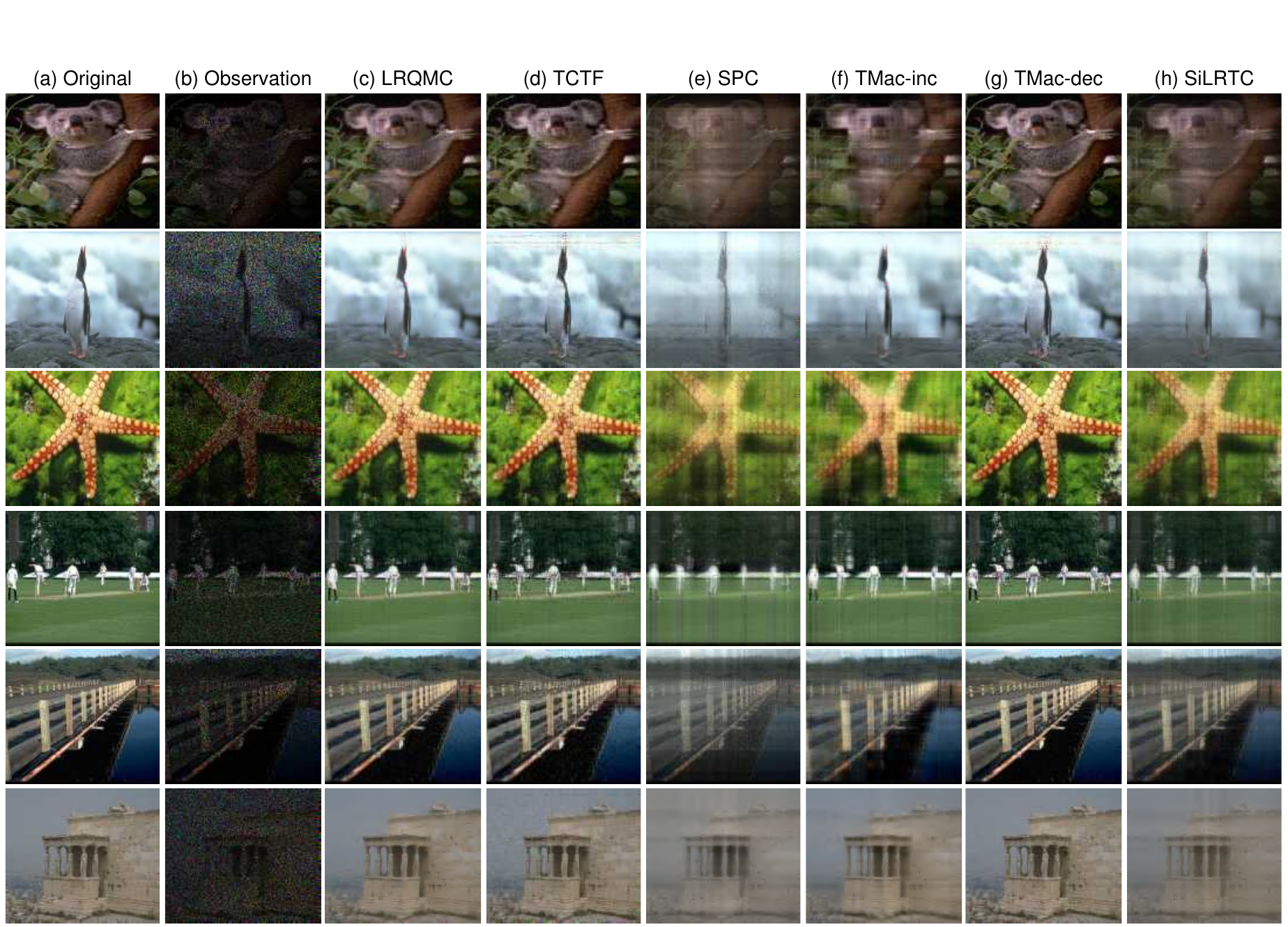}
	\caption{Examples of color image recovery using different algorithms (${\rm{SR}}=0.3$).}
	\label{Examples_of_inpainting}
\end{figure*}

\textbf{Simulation 1:} In this simulation, we use BSD dataset to evaluate our algorithm for color image recovery. We randomly select $56$ color images from this dataset. $6$ examples of selected images are shown in Fig. \ref{Examples_of_inpainting} (a) (from top to bottom, we label them orderly as Image (1), Image (2), Image (3), Image (4), Image (5) and Image (6)) Fig. \ref{Examples_of_inpainting} (b) is the observed image with ${\rm{SR}}=0.3$. Fig. \ref{Examples_of_inpainting} (c)-(h) are the recover results of LRQMC, TCTF, SPC, TMac-inc, TMac-dec and SiLRTC. We see from the Fig. \ref{Examples_of_inpainting} that the color images recovered by LRQMC are visually better than those recovered by the other compared algorithms. TABLE \ref{Examples_of_inpainting_biao} summaries the RSE, PSNR, SSIM, FSIM  values and the running time of all the algorithms on the six testing images displayed in Fig. \ref{Examples_of_inpainting} (a). From the results, one can observe that the overall performance of LRQMC is much better than that of TCTF, SPC TMac-inc and SiLRTC, and also has an advantage over TMac-dec  on most images, except for the Image (5). However, the running time of TMac-dec is much longer than that of LRQMC.

\begin{table*}[htbp]\scriptsize
	\caption{Quantitative quality indexes and running time (seconds) of different algorithms on the six images displayed in Fig. \ref{Examples_of_inpainting} (a)  (${\rm{SR}}=0.3$).}
	\centering
\begin{tabular}{|c|c|c|c|c|c|c|c|}
	\hline
	Images & \diagbox{Indexes}{Algorithms}  & LRQMC  & TCTF  & SPC  & TMac-inc  & TMac-dec  & SiLRTC \\ \toprule
	\hline
	
\multirow{5}{*}{Image (1)}
&RSE &\textbf{-8.519}  &-6.049  &-5.031 &-5.535 &-8.039    &-5.113    \\
\cline{2-8}
&PSNR  &\textbf{29.046}  &23.408  &19.805  &21.893  &28.102    &21.049  \\
\cline{2-8}
&SSIM  &\textbf{0.872}  &0.671  &0.538  &0.631  &0.814    &0.617  \\
\cline{2-8}
&FSIM  &\textbf{0.995}  &0.979  &0.957  &0.962  &0.991    &0.968  \\
\cline{2-8}
&time(s)  &9.578  &9.000  &14.431  &10.571     &50.693  &10.592 \\
	\hline
	
\multirow{5}{*}{Image (2)}	
& RSE &\textbf{-12.858}  &-9.790  &-8.578  &-11.168  &-10.868  &-10.049  \\
\cline{2-8}
&PSNR  &\textbf{28.673}  &24.153  &20.166  &25.293   &24.693    &23.055  \\
\cline{2-8}
&SSIM  &\textbf{0.875}  &0.761  &0.678  &0.808  &0.844   &0.755  \\
\cline{2-8}
&FSIM  &\textbf{0.991}  &0.965  &0.908  &0.953  &0.985    &0.934  \\
\cline{2-8}
&time(s)  &7.421  &16.507  &3.456  &8.900  &61.362   &13.642  \\
\hline
	
\multirow{5}{*}{Image (3)}	
& RSE &\textbf{-9.213}  &-7.108  &-5.100  &-6.107  &-9.138   &-6.186  \\
\cline{2-8}
&PSNR  &\textbf{25.484}  &21.474  &17.177  &21.190  &25.251   &19.348 \\
\cline{2-8}
&SSIM  &\textbf{0.893}  &0.830  &0.718  &0.792  &0.887   & 0.809 \\
\cline{2-8}
&FSIM  &\textbf{0.991}  &0.980  &0.946  &0.956  &\textbf{0.991}   &0.965  \\
\cline{2-8}
&time(s) &16.419  &11.859  &62.164  &8.931  &65.614   &12.765  \\
\hline
	
\multirow{5}{*}{Image (4)}	
& RSE &\textbf{-7.257}  &-4.753  &-5.396  &-6.489  &-5.489   &-5.989  \\
\cline{2-8}
&PSNR  &\textbf{23.350}  &19.504  &19.661  &21.813  &19.814   &20.814  \\
\cline{2-8}
&SSIM  &\textbf{0.791}  &0.678  &0.669  &0.741  &0.786   &0.724  \\
\cline{2-8}
&FSIM  &\textbf{0.988} &0.973  &0.963  &0.967  & 0.982   &0.971 \\
\cline{2-8}
&time(s) &11.863  &9.391  &15.194  &8.168  &64.076   &11.832  \\
\hline
	
\multirow{5}{*}{Image (5)}	
& RSE &-8.097  &-6.409  &-5.141  &-6.143  &\textbf{-8.639 }  &-6.009  \\
\cline{2-8}
&PSNR  &24.480  &21.136  &18.177  &20.176  &\textbf{25.088}   &19.909  \\
\cline{2-8}
&SSIM  &0.706  &0.589  &0.438  &0.564  &\textbf{0.768}   &0.580  \\
\cline{2-8}
&FSIM  &\textbf{0.989}  &0.979  &0.963  &0.967  &\textbf{0.989 }  &0.974  \\
\cline{2-8}
&time(s)  &9.753  &7.048  &7.715  &5.780  &58.974  &12.126   \\
\hline
	
\multirow{5}{*}{Image (6)}	
& RSE  &\textbf{-12.315}  &-9.261  &-8.476  &-10.481  &-12.268    &-9.105  \\
\cline{2-8}
&PSNR  &\textbf{29.613}  &23.425  &21.354  &25.364  &28.941   &22.612  \\
\cline{2-8}
&SSIM  &\textbf{0.901}  &0.724  &0.505  &0.818  &0.895   &0.650  \\
\cline{2-8}
&FSIM  &\textbf{0.989}  &0.973  &0.952  &0.965  &0.986   &0.958  \\
\cline{2-8}
&time(s)  &6.874  &7.645  &13.713  &5.702  &47.972  &15.129   \\
\hline
\end{tabular}
\label{Examples_of_inpainting_biao}
\end{table*}

\begin{figure*}[htbp]
	\centering
	\includegraphics[width=16.5cm,height=16cm]{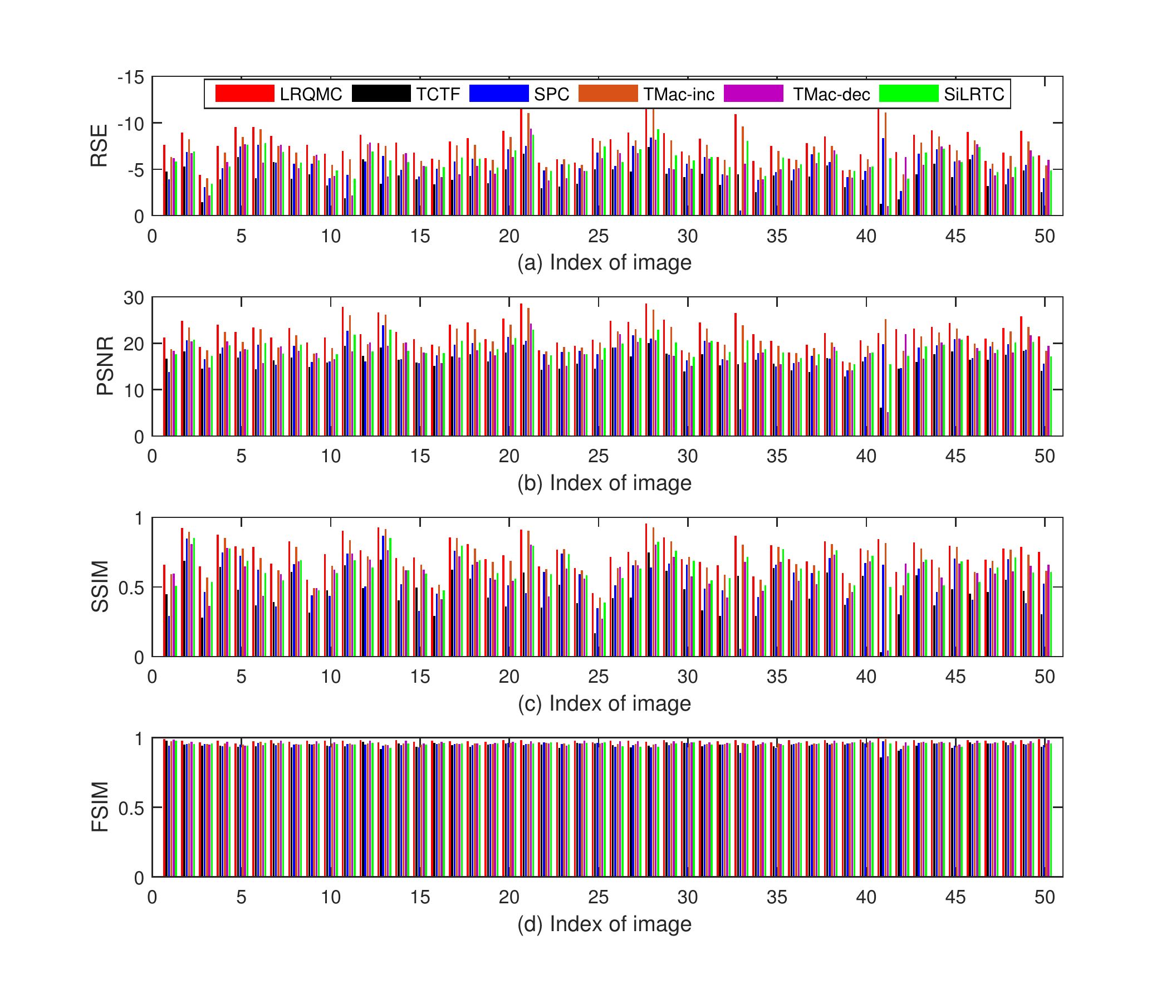}
	\caption{Comparison of RSE, PSNR, SSIM and FSIM results of different algorithms for color image recovery on $50$ BSD images (${\rm{SR}}=0.2$). The figure is viewed better in zoomed PDF.}
	\label{Index4zhu}
\end{figure*}

In Fig. \ref{Index4zhu},  we report the RSE, PSNR, SSIM and FSIM results of different algorithms on the remaining $50$ images. From the results, one can obviously find that our LRQMC algorithm perform better than all the other algorithms in the vast majority of images.

\begin{figure*}[!htb]
	\centering
	\subfigure[]{
		\begin{minipage}{7cm}
			\centering
			\includegraphics[width=6.2cm,height=4cm]{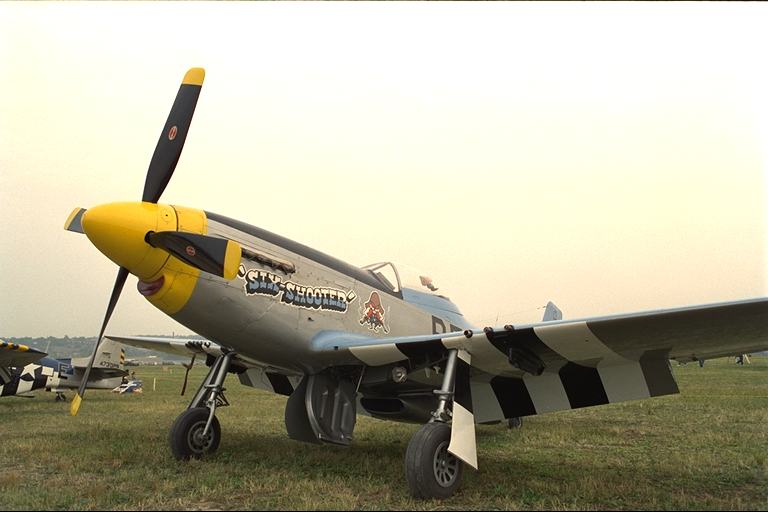}
		\end{minipage}%
	}%
	\subfigure[]{
		\begin{minipage}{7cm}
			\centering
			\includegraphics[width=6.2cm,height=4cm]{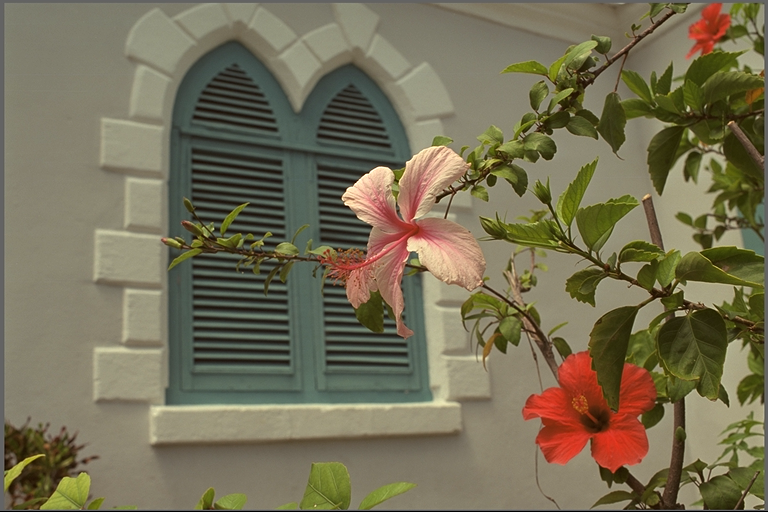}
		\end{minipage}%
	}%
\caption{Randomly selected $2$ color images from Kodak dataset.}
\label{Kodak_tu}
\end{figure*}

\begin{figure*}[htbp]
	\centering
	\includegraphics[width=16.5cm,height=16cm]{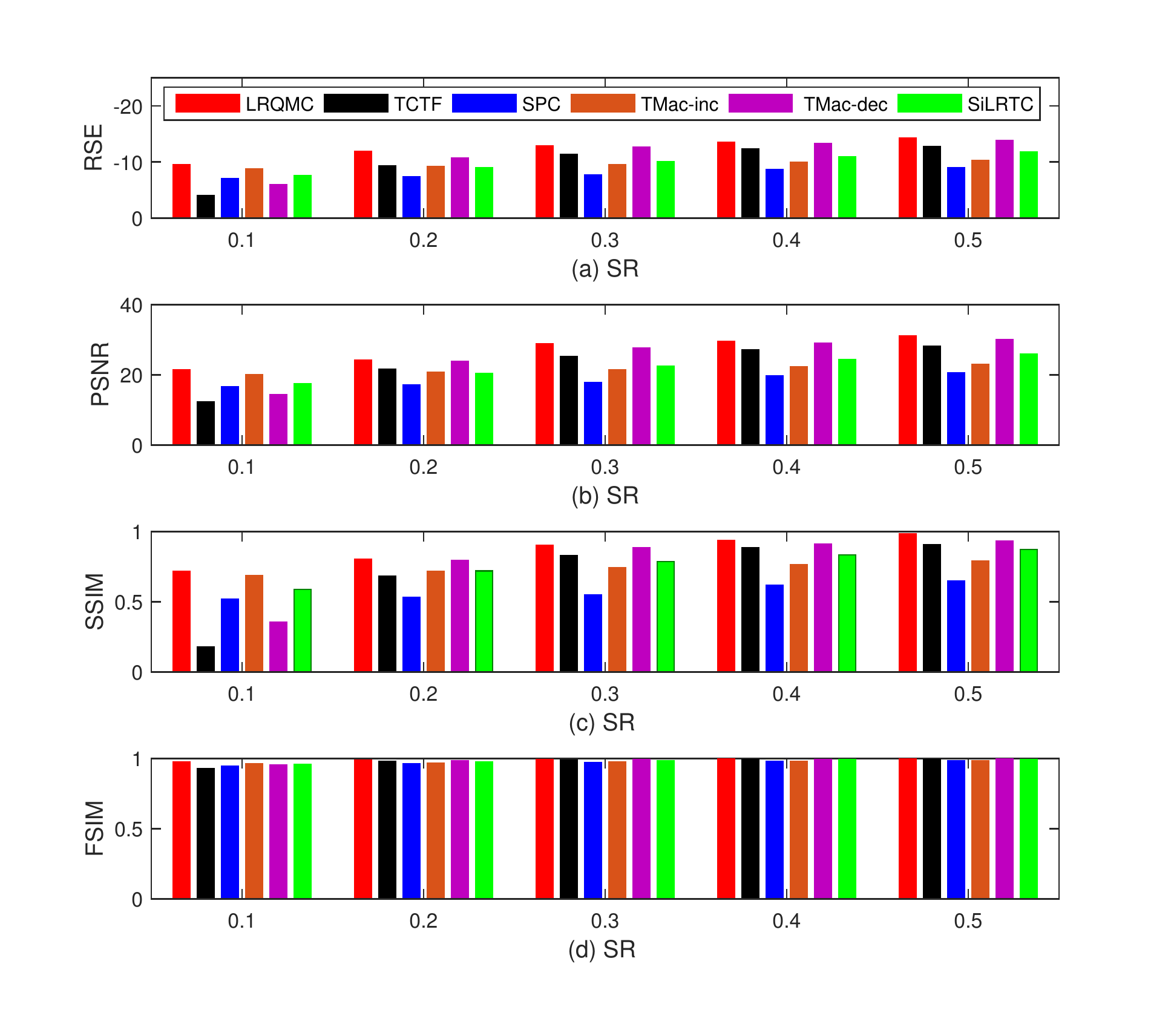}
	\caption{Comparison of RSE, PSNR, SSIM and FSIM results of different algorithms for color image recovery on Fig \ref{Kodak_tu} (a) (${\rm{SR}}=[0.1, 0.2, 0.3, 0.4, 0.5]$). (a) RSE values, (b) PSNR values, (c) SSIM values, (d) FSIM values.}
	\label{Index4SR_1}
\end{figure*}

\textbf{Simulation 2:} In this simulation, we use Kodak dataset to evaluate our algorithm for color image recovery. We randomly select $2$ color images from Kodak dataset as shown in Fig. \ref{Kodak_tu}. Fig. \ref{Index4SR_1} and Fig. \ref{Index4SR_1} respectively show the recovery results of Fig. \ref{Kodak_tu} (a) and Fig. \ref{Kodak_tu} (b) SRs
from $0.1$ to $0.5$ using defferent algorithms in terms of the RSE, PSNR, SSIM and FSIM. LRQMC outperforms all the other algorithms. Furthermore, we can also find that the merit of LRQMC is more obvious in low sampling ratio, \emph{e.g.}, ${\rm{SR}}=0.1$.

\begin{figure*}[htbp]
	\centering
	\includegraphics[width=16.5cm,height=16cm]{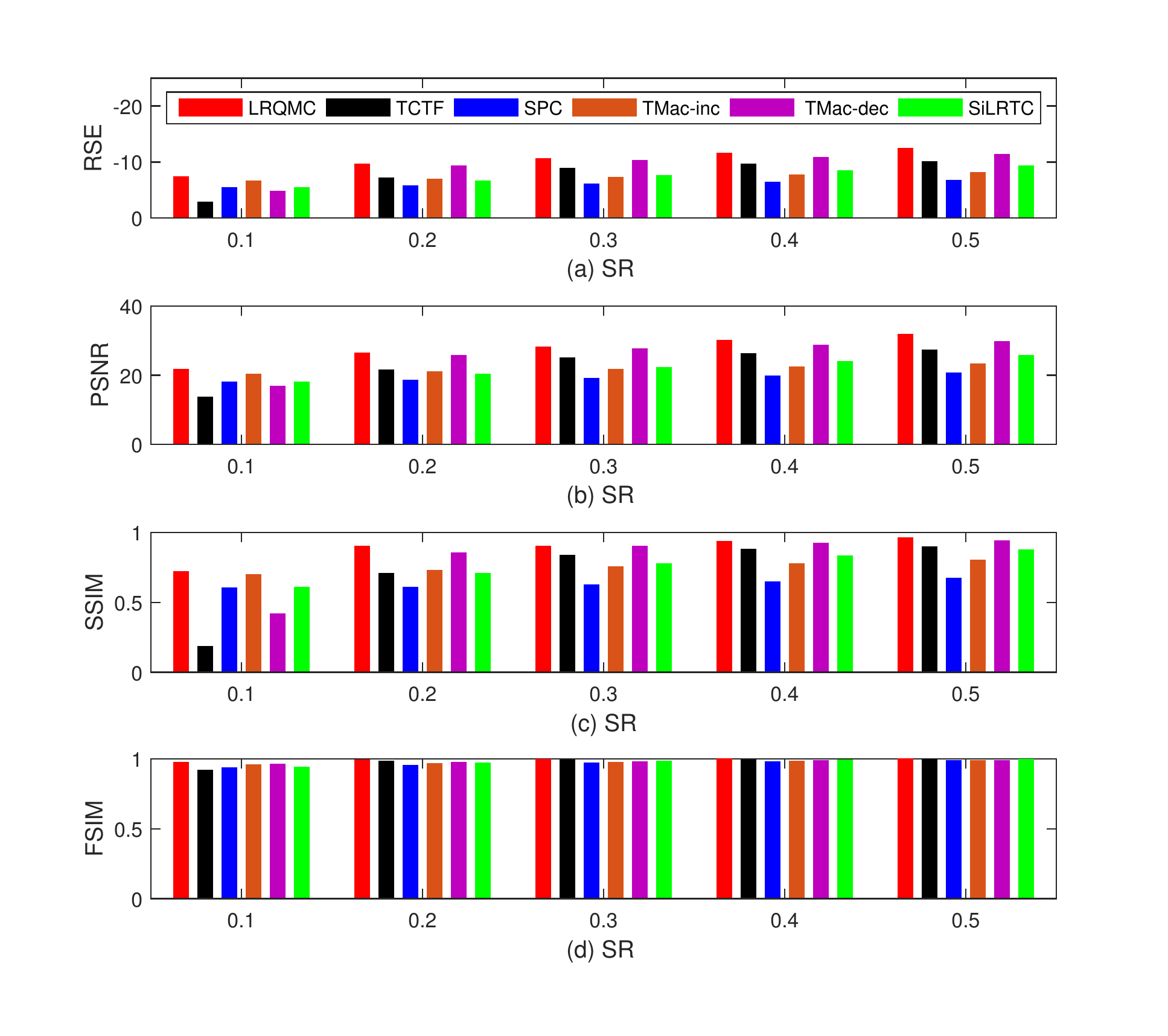}
	\caption{Comparison of RSE, PSNR, SSIM and FSIM results of different algorithms for color image recovery on Fig \ref{Kodak_tu} (b) (${\rm{SR}}=[0.1, 0.2, 0.3, 0.4, 0.5]$). (a) RSE values, (b) PSNR values, (c) SSIM values, (d) FSIM values.}
	\label{Index4SR_2}
\end{figure*}

\section{Conclusion}
We propose a novel low-rank quaternion matrix completion algorithm to recover
missing data of color image. Quaternion representation processes a
colour image holistically as a vector field and handles
the coupling between the color channels naturally, and color information of source image is fully used. We combine low-rank decomposition and nuclear norm (which is replaced by Frobenius norm of the two low-rank factor
quaternion matrices) minimization approaches in our quaternion matrix-based model. Based on the relationship between the quaternion matrix and its equivalent complex matrix, the problem eventually be converted from quaternion number field to complex number field. An alternating minimization method is applied to solve the model, which guarantees convergence of the proposed algorithm. Simulation results on real world color image recovery demonstrate the competitive performance of the proposed algorithm compared to several state-of-the-art tensor-based methods.

\appendices
\section{Proof of the Theorem \ref{theorem2}}
\label{appendices1}
\begin{lemma}
	\label{Q_SVD}
(The SVD of quaternion matrix (QSVD) \cite{10029950538}) Let $\ddot{\mathbf{P}}\in\mathbb{H}^{M\times N}$ be of rank $K$. Then there exist unitary quaternion matrices\footnote{A unitary quaternion matrix $\ddot{\mathbf{B}}\in\mathbb{H}^{N\times N}$ has the following
	property: $\ddot{\mathbf{B}}\ddot{\mathbf{B}}^{H}=\ddot{\mathbf{B}}^{H}\ddot{\mathbf{B}}=\mathbf{I}_{N}$, with $\mathbf{I}_{N}\in\mathbb{R}^{N\times N}$ the identity matrix \cite{DBLP:journals/sigpro/BihanM04}.} $\ddot{\mathbf{A}}\in\mathbb{H}^{M\times M}$ and $\ddot{\mathbf{B}}\in\mathbb{H}^{N\times N}$ such that	
\begin{equation*}
\ddot{\mathbf{A}}\ddot{\mathbf{P}}\ddot{\mathbf{B}}=\left(\begin{array}{cc}
\mathbf{\Sigma}_{K}	&  \mathbf{0}\\
\mathbf{0}	& \mathbf{0}
\end{array}\right),
\end{equation*}
where $\mathbf{\Sigma}_{K}={\rm{diag}}(\sigma_{1},\ldots,\sigma_{K})$ is a real diagonal matrix and has $K$  positive entries $\sigma_{k}, \,(k=1,\ldots,K)$ on its diagonal (\emph{i.e.} positive singular values of $\ddot{\mathbf{P}}$).	
\end{lemma}

According to Lemma \ref{Q_SVD}, we have
\begin{align*}
f(\ddot{\mathbf{A}})f(\ddot{\mathbf{P}})f(\ddot{\mathbf{B}})&=f(\ddot{\mathbf{A}}\ddot{\mathbf{P}}\ddot{\mathbf{B}})\\
&=\left(\begin{array}{cccc}
\mathbf{\Sigma}_{K}& \mathbf{0} &\mathbf{0}  &\mathbf{0}  \\
\mathbf{0}&\mathbf{0}  &\mathbf{0}  &\mathbf{0}  \\
\mathbf{0}&\mathbf{0}  & \mathbf{\Sigma}_{K} & \mathbf{0} \\
\mathbf{0}& \mathbf{0} &\mathbf{0}  & \mathbf{0}
\end{array} \right).
\end{align*}
Hence, from above, it is obvious that ${\rm{rank}}(f(\ddot{\mathbf{P}}))={\rm{rank}}(f(\ddot{\mathbf{A}})f(\ddot{\mathbf{P}})f(\ddot{\mathbf{B}}))=2K$, \emph{i.e.}, ${\rm{rank}}(\ddot{\mathbf{P}})=\frac{1}{2}{\rm{rank}}(f(\ddot{\mathbf{P}}))$ in
Theorem \ref{theorem2} holds.

\section{Proof of the Theorem \ref{proposition2}}
\label{appendices2}
The goal of this appendix is to prove the properties presented
in Theorem \ref{proposition2}.

Proof of (1): According to the QSVD in Lemma \ref{Q_SVD}, there exist unitary quaternion matrices such that
\begin{equation*}
\ddot{\mathbf{X}}=\ddot{\mathbf{A}}\left(\begin{array}{cc}
\mathbf{\Sigma}_{K}	&  \mathbf{0}\\
\mathbf{0}	& \mathbf{0}
\end{array}\right)\ddot{\mathbf{B}}.
\end{equation*}
We let
\begin{equation*}
\ddot{\mathbf{A}}=\left( \begin{array}{cc}
\ddot{\mathbf{A}}_{1}& \ddot{\mathbf{A}}_{2}
\end{array} \right),\;
\ddot{\mathbf{B}}=\left( \begin{array}{c}
\ddot{\mathbf{B}}_{1}\\
\ddot{\mathbf{B}}_{2}
\end{array} \right),
\end{equation*}
where $\ddot{\mathbf{A}}_{1}\in\mathbb{H}^{M\times K}_{K}$, $\ddot{\mathbf{A}}_{2}\in\mathbb{H}^{M\times (M-K)}_{(M-K)}$,
$\ddot{\mathbf{B}}_{1}\in\mathbb{H}^{K\times N}_{K}$, $\ddot{\mathbf{B}}_{2}\in\mathbb{H}^{(N-K)\times N}_{(N-K)}$. Then, we have
\begin{align*}
\ddot{\mathbf{X}}&=\left( \begin{array}{cc}
\ddot{\mathbf{A}}_{1}& \ddot{\mathbf{A}}_{2}
\end{array} \right)
\left(\begin{array}{cc}
\mathbf{\Sigma}_{K}	&  \mathbf{0}\\
\mathbf{0}	& \mathbf{0}
\end{array}\right)
\left( \begin{array}{c}
\ddot{\mathbf{B}}_{1}\\
\ddot{\mathbf{B}}_{2}
\end{array} \right) \\
&=\ddot{\mathbf{A}}_{1}\mathbf{\Sigma}_{K}\ddot{\mathbf{B}}_{1}\\
&=\ddot{\mathbf{U}}\ddot{\mathbf{V}},
\end{align*}
where $\ddot{\mathbf{U}}=\ddot{\mathbf{A}}_{1}\in\mathbb{H}^{M\times K}_{K}$, $\ddot{\mathbf{V}}=\mathbf{\Sigma}_{K}\ddot{\mathbf{B}}_{1}\in\mathbb{H}^{K\times N}_{K}$.

Proof of (2): Recall that if $\mathbf{P}\in \mathbb{C}^{M\times N}$ and $\mathbf{Q}\in \mathbb{C}^{N\times M}$ are two matrices, then we have ${\rm{rank}}(\mathbf{P}\mathbf{Q})\leq {\rm{min}}({\rm{rank}}(\mathbf{P}),{\rm{rank}}(\mathbf{Q}))$. Thus, we immediately have
\begin{align*}
{\rm{rank}}(\ddot{\mathbf{P}}\ddot{\mathbf{Q}})&=\frac{1}{2}{\rm{rank}}(f(\ddot{\mathbf{P}}\ddot{\mathbf{Q}}))\\
&=\frac{1}{2}{\rm{rank}}(f(\ddot{\mathbf{P}})f(\ddot{\mathbf{Q}}))\\
&\leq \frac{1}{2}{\rm{min}}({\rm{rank}}(f(\ddot{\mathbf{P}}), {\rm{rank}}(f(\ddot{\mathbf{Q}}))\\
&={\rm{min}}({\rm{rank}}(\ddot{\mathbf{P}}),  {\rm{rank}}(\ddot{\mathbf{Q}})).
\end{align*}

For property (3), actually, similar result and proof can be found in \cite{DBLP:journals/siamrev/RechtFP10} (Lemma $5.1$) and \cite{ DBLP:journals/jstsp/YangPCO18}.

\section*{Acknowledgment}
This work was supported by The Science and Technology Development Fund, Macau SAR (File no. FDCT/085/2018/A2).




%

\bibliographystyle{elsarticle-num}
\bibliography{bare_jrnl}

%

%
%
%
%
%
%
%
%
%
\end{document}